\documentclass[12pt,amsfonts,amssymb]{article}
\usepackage[english]{babel}

\ifx\relax\mathfrak \@xp\@gobbletwo 
  \else \let\mathfrak\relax \fi 
\RequirePackage{amsfonts}\relax 
\RequirePackage{amsmath}\relax 

\RequirePackage{oldlfont}
\def\carre{{$\Box$}}

\makeatletter 
\def\@ypreuve[#1]{\@preuve{ #1}} 
\def\@preuve#1{\begin{trivlist}\item[]{\em Proof#1: }} 
\newenvironment{preuve}{\@ifnextchar[{\@ypreuve}{\@preuve{}}}{\hfill\carre\end{trivlist}} 

\newtheorem{theorem}{Theorem} 
 
\newtheorem{lemma}[theorem]{Lemma}

\def\re{{\mathrm{Re}\,}} 
\def\im{{\mathrm{Im}\,}} 
 
\def\dom{{\mathrm{dom}}}

\begin{document}
\begin{center}
{\LARGE\bf Note on a product formula \\ for unitary groups}
\bigskip
\bigskip

by Vincent CACHIA\footnote{Department of Theoretical 
Physics, quai Ernest-Ansermet 24, CH-1211 GENEVA 4 SWITZERLAND. 
email: Vincent.Cachia@physics.unige.ch}
\medskip
\date{\today}
\bigskip
\end{center}

\begin{quote}
{\small
\noindent{\bf Abstract:} 
For any nonnegative self-adjoint operators $A$ and $B$ in a separable Hilbert space,
we show that the Trotter-type formula
$[(e^{i2tA/n}+e^{i2tB/n})/2]^n$ converges strongly in $\overline{\dom(A^{1/2})\cap\dom(B^{1/2})}$
for some subsequence and for almost every $t\in\mathbb R$. This result extends to the degenerate case 
and to Kato-functions following the method of Kato \cite{Kato78}.}
\end{quote}

\bigskip
\bigskip

In a famous paper \cite{Kato78}, T. Kato proved that for any nonnegative self-adjoint operators
$A$ and $B$ in a Hilbert space $\cal H$, the Trotter product formula 
$(e^{-tA/n}e^{-tB/n})^n$ converges strongly to the (degenerate)
semigroup generated by the form-sum $A\dot{+}B$ for any $t$ with $\re t>0$. He also extended
the result to a class of so-called Kato-functions, and to degenerate semigroups.
However the convergence on the boundary $i\mathbb R$ remains an unclear problem 
in this generality \cite{Chernoff,Exner}. For Kato-functions $f$ such that $\im f\leq 0$ (for example
$f(s)=(1+is)^{-1}$), Lapidus found such an extension \cite{Lapidus}. For the case of the Trotter
product formula with projector $(e^{itA/n}P)^n$, Exner and Ichinose obtained recently a
interesting result \cite{ExIc}.

\section{Statement of the result}

Since this note is closely related to Kato's paper \cite{Kato78}, it is convenient to use
similar notations. $A$ and $B$ denote nonnegative self-adjoint operators defined in closed subspaces
$M_A$ and $M_B$ of a separable Hilbert space $\cal H$, and $P_A$, $P_B$ denote the orthogonal
projections on $M_A$ and $M_B$. Let ${\cal D}'=\dom(A^{1/2})\cap\dom(B^{1/2})$, let $\cal H'$
be the closure of $\cal D'$, and let $P'$ be the orthogonal projector on $\cal H'$.
The form-sum $C=A\dot{+}B$ is defined as the self-adjoint operator
in $\cal H'$ associated with the nonnegative, closed quadratic form 
$u\mapsto \|A^{1/2}u\|^2 + \|B^{1/2}u\|^2$, $u\in{\cal D'}$. 
We consider Trotter-type product formulae $F(t/n)^n$ based on the arithmetic mean
\begin{equation}
F(t) = \frac{f(2tA)P_A + g(2tB)P_B}{2}.
\end{equation}
The Kato-functions $f$ and $g$ are assumed here to be bounded, holomorphic functions
in $\{t\in{\mathbb C}: \re t> 0\}$ with:
\begin{equation}
 |f(t)|\leq 1,\ f(0)=1,\ f'(+0)=\lim_{t\rightarrow 0,\re t>0}\frac{f(t)-1}{t} = -1,
\end{equation}
and $0\leq f(s)\leq 1$ if $s>0$, and the same conditions for $g$.
%may be either $s\mapsto e^{-s}$ or $s\mapsto (1+s)^{-1}$,
By the functional calculus for normal operators, $F(t)$ is well
defined for $\re t\geq 0$ and bounded by $1$.

\begin{theorem}
Let $\cal H$ be a separable Hilbert space.
Let $A$, $M_A$, $P_A$, $B$, $M_B$, $P_B$, $C$, $\cal H'$, $P'$, $f$, $g$, and $F$ be as defined
above.
For any $u\in\cal H$ one has
\begin{equation}\label{th1.2}
\lim_{n\rightarrow\infty} \int_{-\infty}^{+\infty} \phi(t)F(it/n)^n udt = \int_{-\infty}^{+\infty}
\phi(t)e^{-itC}P'udt,\ \phi\in L^1(\mathbb R).
\end{equation}
Moreover there exists a set $L\subset\mathbb R$ with zero Lebesgue measure and an increasing
function $\varphi:\mathbb N\rightarrow N$, such that:
\begin{equation}\label{th1.1}
\forall u\in{\cal H'},\quad 
F(it/\varphi(n))^{\varphi(n)} u \longrightarrow e^{-itC}P'u,\ t\in{\mathbb R}\setminus L,
\end{equation}
whereas for $u\in M_A^\perp+M_B^\perp$, $\lim_{n\rightarrow\infty} F(it/n)^n u=0$,
for each $t\in\mathbb R$.
\end{theorem}

One sees that the strong convergence, valid in the open right half-plane, cannot extend exactly
to the boundary $i\mathbb R$, as already remarked in \cite{Exner, ExIc}: the strong
convergence on the boundary is restricted to the subspace ${\cal H'}+M_A^\perp+M_B^\perp$.
However the weaker convergence (\ref{th1.2}) was already observed \cite{Chernoff,Feldman}.

\section{Proof}

Let us consider for $\re t\geq 0$ and $\tau>0$:
\begin{equation}
S_{t,\tau}=\tau^{-1}(I-F(t\tau))
\end{equation}
which is a holomorphic operator-valued function of $t$ in the open right half-plane.

The main step of the proof is to show that the strong convergence:
\begin{equation}
s-\lim_{\tau\rightarrow 0} (I+S_{t,\tau})^{-1} = (I+tC)^{-1}P'
\end{equation}
holds for $\re t>0$, and remains true for almost all $t\in i\mathbb R$ and on some subsequence.
This will give the announced result (\ref{th1.1}) for $u\in{\cal H'}$ by Chernoff's lemma
(see below). The convergence on the boundary is obtained by a useful result of Feldman
\cite[Th 5.1]{Feldman}, which we state here in a slightly more general form:

\begin{lemma}\label{feldman}
Let $\{\Psi_\tau : 0<\tau<1\}$ be a uniformly bounded family of bounded holomorphic
$\cal H$-valued function defined in the open right half-plane.
Suppose that 
$\Psi_\tau(z)\stackrel{\tau\rightarrow 0}{\mathop{\longrightarrow}}
\Psi(z)$, for each $z$ with $\re z>0$. 
%Then $\sigma(L^\infty,L^1)-\lim_{\tau\rightarrow 0} \Psi(i\cdot)=\Psi(i\cdot)$, which means: 
Then for each $v\in L^1(\mathbb R,\cal H)$
\begin{equation}\label{vectfeldman}
\int_{\mathbb R} (v(t),\Psi_\tau(it))dt
\stackrel{\tau\rightarrow 0}{\mathop{\longrightarrow}}
\int_{\mathbb R} (v(t),\Psi(it))dt.
\end{equation}
Moreover for each (numerical) $\phi\in L^1(\mathbb R)$:
\begin{equation}\label{numfeldman}
s-\lim_{\tau\rightarrow 0} \int_{\mathbb R} \phi(t)\Psi_\tau(it)dt
= \int_{\mathbb R} \phi(t)\Psi(it)dt.
\end{equation}
\end{lemma}

\begin{preuve}
The two results (\ref{vectfeldman}) and(\ref{numfeldman}) follow from very similar arguments, 
so we present only the first one.
The bounded holomorphic functions $\Psi_\tau$ have boundary values for almost every $is\in i\mathbb R$,
and for any $t>0$ and $s\in\mathbb R$,
$$
\Psi_\tau(t+is) = \int_{-\infty}^{+\infty}\frac{t/\pi}{t^2+(s-s')^2}\Psi_\tau(is')ds' 
= P_t\ast\Psi_\tau(i\cdot).
$$
The kernel $P_t$ is in fact an approximate identity: $P_t\ast\phi
\stackrel{L^1}{\mathop{\longrightarrow}}\phi$ as $t\rightarrow 0$. Then we have:
\begin{eqnarray*}
\int_{\mathbb R}(v(s),[\Psi_\tau(is)-\Psi(is)])ds & = & 
\int_{\mathbb R}(v(s),[\Psi_\tau(t+is)-\Psi(t+is)])ds + \\
& & + \int_{\mathbb R} ([v(s)-(P_t\ast v)(s)],[\Psi_\tau(is)-\Psi(is)])ds,
\end{eqnarray*}
where we have used:
$$\int_{\mathbb R} (v(s),[P_t\ast h](s))ds = \int_{\mathbb R} ([P_t\ast v](s),h(s))ds.$$
This leads to
\begin{eqnarray}
\left|\int_{\mathbb R}(v(s),[\Psi_\tau(is)-\Psi(is)])ds\right| & \leq & 
\int_{\mathbb R}\|v(s)\| \|\Psi_\tau(t+is)-\Psi(t+is)\|ds \label{intfeld} \\
& &\hspace{-2cm} + \int_{\mathbb R}\|v(s)-(P_t\ast v)(s)\| \|\Psi_\tau(is)-\Psi(is)\|ds.\label{last}
\end{eqnarray}
The last term (\ref{last}) can be made arbitrary small by choosing $t$ sufficiently
small, and for any $t$ the integral in the right hand side of (\ref{intfeld}) tends to $0$
as $\tau\rightarrow 0$ by Lebesgue's theorem.

\end{preuve}

\noindent{\bf Remark:} It can also be shown that (\ref{numfeldman}) implies (\ref{vectfeldman}).
\bigskip

It is convenient to introduce the following bounded accretive operators, for $\re t\geq 0$
and $\tau>0$:
\begin{equation}
A_{t,\tau} = \tau^{-1}[I-f(t\tau A)P_A],\mbox{ and } B_{t,\tau} =\tau^{-1}[I-g(t\tau B)P_B].
\end{equation}

\begin{lemma}\label{lem1}
For any $t$ with $\re t>0$, $s-\lim_{\tau\rightarrow 0} (I+S_{t,\tau})^{-1} = (I+tC)^{-1}P'$.
Moreover for any $v\in L^1({\mathbb R},{\cal H})$, $u\in\cal H$ and $t\in\mathbb R$, one has
\begin{equation}\label{wconv1}
\lim_{\tau\rightarrow 0} \int_{\mathbb R} (v(t),(I+S_{it,\tau})^{-1}u) dt =
\int_{\mathbb R} (v(t),(I+itC)^{-1}P'u)dt.
\end{equation}
\end{lemma}

\begin{preuve}
Since $S_{t,\tau} = A_{t,2\tau}+B_{t,2\tau}$, the strong convergence of $(I+S_{t,\tau})^{-1}$
for $t>0$ follows from
\cite[Lem. 2.2 and 2.3]{Kato78}. Then it extends to the open right half-plane by the theorem
of Vitali: for any $\tau>0$, $(I+S_{t,\tau})^{-1}$ is a holomorphic function of $t$, and is
bounded by $1$.

The convergence on the boundary (\ref{wconv1}) follows from Lemma \ref{feldman}.
\end{preuve}

For any fixed $u\in{\cal H}$ and $t\in\mathbb R$ we set $w_{t,\tau}=(I+S_{it,\tau})^{-1}u$, $\tau>0$.
Then one finds
\begin{equation}\label{eq1}
(u,w_{t,\tau}) = \|w_{t,\tau}\|^2 + (A_{it,2\tau} w_{t,\tau},w_{t,\tau}) +
(B_{it,2\tau} w_{t,\tau},w_{t,\tau})
\end{equation}
with $\re (A_{it,2\tau}w_{t,\tau},w_{t,\tau})\geq 0$ and
 $\re (B_{it,2\tau}w_{t,\tau},w_{t,\tau})\geq 0$. Therefore
\begin{equation}\label{estiw}
\|w_{t,\tau}\|^2 \leq \re(u,w_{t,\tau}) \leq |(u,w_{t,\tau})| \leq \|u\| \|w_{t,\tau}\|
\end{equation}
and thus $\|w_{t,\tau}\| \leq \|u\|$, $\tau>0$.

\begin{lemma}\label{lem4}
Let $\alpha_n$ be any sequence of positive numbers with limit zero.
There exists a set $L\subset\mathbb R$ of zero Lebesgue measure and a subsequence $\tau_n$ of
$\alpha_n$, such that for each  $t\in{\mathbb R}\setminus L$,
$s-\lim_{n\rightarrow\infty}(I+S_{it,\tau_n})^{-1}=(I+itC)^{-1}P'$.
\end{lemma}

\begin{preuve}
It follows from Lemma \ref{lem1} that:
$\int_{\mathbb R}(1+t^2)^{-1} (u,w_{t,\tau})dt \rightarrow
\int_{\mathbb R} (1+t^2)^{-1}(u,w_t)dt$ with $w_t = (I+itC)^{-1}P'u$.
Thus the same is true for the real part, and we have by (\ref{eq1}):
\begin{equation}
\re(u,w_{t,\tau}) = \|w_{t,\tau}\|^2+\|(\re A_{it,2\tau})^{1/2}w_{t,\tau}\|^2
+\|(\re B_{it,2\tau})^{1/2}w_{t,\tau}\|^2 \nonumber
\end{equation}
We observe that $\re(u,w_t) = \re ((I+itC)w_t,w_t) = \|w_t\|^2$, and that
$$
\int_{\mathbb R}\re(w_{t,\tau},w_t)\frac{dt}{1+t^2}
 \stackrel{\tau\rightarrow 0}{\mathop{\longrightarrow}}
\int_{\mathbb R} \|w_t\|^2 \frac{dt}{1+t^2}.
$$
Then one finds
$$
\int_{\mathbb R} \left(\|w_{t,\tau}-w_t\|^2 + \|(\re A_{it,2\tau})^{1/2}w_{t,\tau}\|^2
+\|(\re B_{it,2\tau})^{1/2}w_{t,\tau}\|^2\right)\frac{dt}{1+t^2}
\stackrel{\tau\rightarrow 0}{\mathop{\longrightarrow}} 0,
$$
in particular $\int_{\mathbb R} \|w_{t,\tau}-w_t\|^2 (1+t^2)^{-1}dt
\stackrel{\tau\rightarrow 0}{\mathop{\longrightarrow}} 0$. This means that the functions
$t\mapsto\|w_{t,\tau}-w_t\|$ converge to $0$ in $L^2(\mathbb R,\mu)$ as $\tau\rightarrow 0$,
with the finite measure $d\mu=(1+t^2)^{-1}dt$.
Let $(e_m)_{m\in\mathbb N}$ be a basis of the separable Hilbert space $\cal H$.
For $u=e_1$ the above $L^2$-convergence implies that there exists $L_1\subset\mathbb R$
with $\mu(L_1)=0$ and
some increasing function $\varphi_1:\mathbb N\rightarrow N$ such that
$(I+S_{it,\alpha_{\varphi_1(n)}})^{-1}e_1 \rightarrow (I+itC)^{-1}P'e_1$ as $n\rightarrow\infty$,
for any $t\in{\mathbb R}\setminus L_1$.
Then for $u=e_2$ there exists $L_2\subset\mathbb R$ with $\mu(L_2)=0$ and an increasing
function $\varphi_2:\mathbb N\rightarrow\mathbb N$, such that
$(I+S_{it,\alpha_{\varphi_1\circ\varphi_2(n)}})^{-1}e_2 \rightarrow (I+itC)^{-1}P'e_2$ as 
$n\rightarrow\infty$,
for any $t\in{\mathbb R}\setminus L_2$, and so on for each $m\in\mathbb N$.
Finally by the diagonal procedure we consider the sequence 
$\tau_n=\alpha_{\varphi_1\circ\cdots\circ\varphi_n(n)}$
and find that convergence holds for each vector $e_m$ of the basis and for each 
$t\in{\mathbb R}\setminus L$, where $L=\cup_{m\in\mathbb N} L_m$. We have $\mu(L)=0$ and thus
$L$ has also zero Lebesgue measure.
Since the operators $(I+S_{it,\tau})^{-1}$ are uniformly bounded, this implies the strong
convergence for any vector $u\in{\cal H}$.

\end{preuve}

\begin{preuve}[of the theorem]
We consider $Z_{t,n}=(n/t)[F(it/n)-I] = -t^{-1}S_{it,1/n}$ and $\alpha_n=1/n$.
Let $L$ be as in Lemma \ref{lem4} and let $t\in{\mathbb R}\setminus L$, $t\neq 0$. 
By \cite[Th. 3.17]{Davies} and Lemma \ref{lem4} one obtains for some increasing function 
$\varphi:\mathbb N\rightarrow N$:
\begin{equation}
\lim_{n\rightarrow\infty}e^{sZ_{t,\varphi(n)}}u = e^{-isC}P'u,\quad u\in{\cal H}',\ s\in\mathbb R.
\end{equation}
By Chernoff's Lemma \cite[Lem. 3.27 and 3.29]{Davies}, one has
\begin{equation}
\lim_{n\rightarrow\infty} \|F(it/\varphi(n))^{\varphi(n)} u - e^{\varphi(n)(F(it/\varphi(n))-I)}u\| = 0,
\quad u\in{\cal H}'.
\end{equation}
Thus we obtain the convergence (\ref{th1.1}) in $\cal H'$. If $u\in M_A^\perp+M_B^\perp$, the convergence
to $0$ is clear because $F(it/n)^n$ reduces to $f(2it/n)^n P_A/2^n$ or $g(2it/n)^n P_B/2^n$.
The $\sigma(L^\infty,L^1)$ convergence (\ref{th1.2}) follows by using Kato's result for $\re t>0$ and
Lemma \ref{feldman}.

\end{preuve}

{\bf Remark:} the subsequence appearing in the theorem makes the result somewhat unsatisfactory. In fact
this restriction is not necessary if we assume that the functions $t\mapsto (I+S_{it,\tau})^{-1}u$ are
equicontinuous with respect to $\tau$, at some point $t_0\neq 0$. 
In this case Lemma \ref{lem4} can be improved in the following way:
\begin{equation}\label{newlem4}
s-\lim_{\tau\rightarrow 0}(I+S_{it_0,\tau})^{-1} = (I+it_0 C)^{-1}P'.
\end{equation}
For the proof, let us consider an approximate identity $\rho_n:\mathbb R\mapsto R_+$. By Lemma~\ref{lem1}
one has $\lim_{\tau\rightarrow 0} [\rho_n\ast(u,w_{\cdot,\tau})](t_0)=[\rho_n\ast(u,w_{\cdot})](t_0)$
for each $n=1,2,\dots$, and by the equicontinuity of the functions $t\mapsto(u,w_{t,\tau})$ at $t_0$,
$\lim_{n\rightarrow\infty}[\rho_n\ast(u,w_{\cdot,\tau})](t_0)=(u,w_{t_0,\tau})$ uniformly in $\tau$.
Then in the proof of the theorem we consider $Z_{t,n}=-t_0^{-1}S_{it_0,t/nt_0}$ for any $t\in\mathbb R$.
 By (\ref{newlem4})
one has $s-\lim_{n\rightarrow\infty}(t_0^{-1}-Z_{t,n})^{-1}=(t_0^{-1}+iC)^{-1}P'$,
which leads to the result of the theorem without subsequence
(the exceptional set $L$ has also disappeared). Concerning the equicontinuity condition, we recall
that in our first result the function $t\mapsto F(it)$ ($t\in\mathbb R$) is not necessarily continuous
(whereas it is continuous for example if $f$ and $g$ are the exponential function).

\end{document}